\documentclass[10pt,a4paper]{article}
\usepackage{amssymb}
\newtheorem{thm}{Theorem}
\newtheorem{defn}{Definition}
\newtheorem{cor}{Corollary}
\newtheorem{lem}{Lemma}

\newtheorem{prop}{Proposition}
\newtheorem{ex}{Example}
\newtheorem{rk}{Remark}

\usepackage{color}
\usepackage{colordvi}
\usepackage{amsmath}
\usepackage{amssymb}

\def \endsquare{ \hfill $\sqcup \!\!\!\! \sqcap$ \par}
\def\cutoffint{-\hskip -10pt\int}
\def \C{{\! \rm \
I \!\!\!C}}
\def \R {{\! \rm \ I \!R}}

\def \N {{\! \rm \ I \!N}}

\def \Z {{\! \rm Z\! \!Z}}
 \def\cutoffint{-\hskip -10pt\int}

\def \e {{\epsilon}}

\def \Ci {{C^\infty}}

\def \Cl {{C\ell}}

\def\res{{\rm res}}

\long\def\ignore#1{}

\begin{document}

\title{\bf  Stokes' formulae on classical symbol valued forms and applications}

\author{Yoshiaki MAEDA, Dominique MANCHON, Sylvie
PAYCHA }

\maketitle
\section*{Abstract}
The  Wodzicki residue and the cut-off integral extend to classical symbol
valued forms.  We show that  they obey a Stokes' type property and that the extended Wodzicki residue can be interpreted as a complex residue  like the ordinary one.  \\
 In the case of cut-off integrals, Stokes' property  (i.e. vanishing on exact
 forms) only  holds for non 
integer order symbol valued forms and leads to an integration by parts formula and translation invariance for cut-off integrals on  non integer order classical symbols. \\ 
The extended Wodzicki residue yields an even residue cycle on classical
symbols and an odd cochain (the {\sl cosphere cochain\/}) which measures an obstruction to Stokes'
property of the cut-off integral on integer order symbol-valued forms.
\section*{R\'esum\'e}
Le r\'esidu de Wodzicki et l'int\'egrale r\'egularis\'ee par
troncature s'\'etendent aux formes \`a coefficients symboles classiques. Nous
montrons que 
que l'un et l'autre poss\`edent une propri\'et\'e de Stokes et que  le r\'esidu
de Wodzicki des formes s'interpr\`ete comme un r\'esidu complexe, de la m\^eme
mani\`ere que le r\'esidu de Wodzicki ordinaire. \\
Dans le cas de l'int\'egrale r\'egularis\'ee par troncature, la propri\'et\'e
de Stokes (i.e. l'annulation sur les formes exactes) n'est v\'erifi\'ee que
pour les formes d'ordre non-entier. Elle implique une formule d'int\'egration
par parties et une invariance par translation pour l'int\'egrale
r\'egularis\'ee des symboles d'ordre non-entier.\\
Le r\'esidu de Wodzicki \'etendu induit quant \`a lui un cycle de dimension
paire sur l'alg\`ebre
engendr\'ee par les symboles classiques, ainsi qu'une cocha\^\i ne de degr\'e
un de moins (la {\sl cocha\^\i ne cosph\`ere\/}) qui mesure l'obstruction \`a
la propri\'et\'e de Stokes pour les formes d'ordre entier.

\section*{Introduction}
We discuss generalisations of Stokes' property $\int_{U} d\alpha=0$ for ordinary integrals
 of forms $\alpha$ with compact support (or tending to zero rapidly enough at infinity) in an open subset $U$ of $\R^n$  to regularised integrals of classical symbol valued forms on an open subset of $\R^n$. Although consequences of such a formula such as  integration by parts  and translation invariance for regularised integrals are commonly used in the physics literature to compute Feynman graphs, the only explicit reference we could find in the literature to Stokes' formula for regularised integrals is in \cite{E}. Etingof considers dimensional regularisation which he applies to a class of functions relevant for physics, namely  functions of Feynman type,  proving Stokes' formula for corresponding regularised integrals of top degree forms.\\ 
Here, we consider general regularisation procedures and all     classical symbol valued forms, proving Stokes' formula with pseudodifferential theoretic tools; an essential obstacle to Stokes' formula turns out to be the Wodzicki residue extended to forms, to which we devote a large part of the paper. \\ \\
The  Wodzicki residue extended to classical symbol valued forms is the topic
of the first part of the paper. It satisfies Stokes' property and therefore
defines  a $2n$-cycle on the algebra of classical symbols with compact support
on an open subset $U\subset\R^n$ equipped with the left product of symbols
$\ast$ (Theorem \ref{thm:residuecycle}). Its associated  residue character
is a cyclic $\star$-Hochschild cocycle: 
$$(\sigma_0, \cdots,\sigma_{2n})\mapsto {\rm res}\left( \sigma_0\ast d\, \sigma_1\wedge_\ast \cdots \wedge_\ast d\,\sigma_{2n}\right)$$
where  res is the extended residue and where  $\wedge_\ast$ is the product on
the graded differential algebra of classical symbol valued forms induced by
the left product on symbols and $d$ the exterior differentiation  on
$T^*U$.\\ 
On classical pseudodifferential operators of order $0$,
  this $\star$-Hochschild cocycle  reduces to a cyclic Hochschild cocycle for
  the ordinary product for we have:  $${\rm res}\left( \sigma_0\ast d\, \sigma_1\wedge_\ast \cdots \wedge_\ast d\,\sigma_{2n}\right)={\rm res}\left( \sigma_0\, d\, \sigma_1\wedge \cdots \wedge d\,\sigma_{2n}\right).$$
 It
  coincides up to a multiplicative constant with the analog in the context of classical symbols of the  antisymmetrised  $2n$-cocycle
introduced in \cite{CFS} and further investigated in \cite{H} in the context of
star-deformed algebras:
  $${\rm res}\left( \sigma_0\ast d\, \sigma_1\wedge_\ast \cdots \wedge_\ast
    d\,\sigma_{2n}\right)= \frac{(-i)^n}{n!}\,A\, {\rm res} \left( \sigma_0 \star
    \theta(\sigma_1, \sigma_2) \star \cdots \star
    \theta(\sigma_{2n-1}, \sigma_{2n})\right)$$
where we have set  $\theta(\sigma_i, \sigma_j)= \sigma_1\star \sigma_j-\sigma_i
\cdot \sigma_j$ as in \cite{H}, \cite{CFS}. Here $A$ is the antisymmetrisation over all but the first
variable.  
\\ \\ The second part of the paper is devoted to cut-off  integrals which we   also extend to  classical symbol valued forms. We show they obey Stokes'
property when restricted to non integer order symbols with compact support (see Theorem \ref{thm:Stokescutoff}). As a result, we get an integration by parts formula for cut-off integrals on non integer order symbols  and show translation invariance for cut-off integrals on non integer order symbols.\\
Stokes' property does not hold anymore on  integer order symbol valued forms  with compact support; as a result, one does not expect to define a cycle on the algebra of classical symbols using cut-off integrals.  Rather, we express  the obstruction to the cyclicity of a $2n$-cochain defined in terms of  cut-off integrals of symbols $\cutoffint_{T^*U}$: 
$$(\sigma_0, \cdots,\sigma_{2n})\mapsto \cutoffint_{T^*U}\left(\sigma_0\ast d\, \sigma_1\wedge_\ast \cdots \wedge_\ast d\,\sigma_{2n}\right)_0$$
where the subscript $0$ stands for the $0$-order part of the symbol valued
form, in terms of the {\sl cosphere $2n-1$-cochain\/} defined in a similar way to the residue character (Proposition \ref{prop:Bcutoff}). 
\\ \\
 Finally, in the third part of the paper, we show  that  the relation between complex residues and the Wodzicki residue extends to symbols valued forms (Theorem \ref{thm:KVforms}):
$${\rm Res}_{z=z_0} \cutoffint \omega(z)=-\frac{1}{\alpha^\prime(z_0)}{\rm res}(\omega(z_0)),$$
where $\omega(z)$ is a holomorphic family of classical symbol valued forms of order $\alpha(z)$ and ${\rm res}(\omega(z_0))$ the Wodzicki residue of $\omega(z_0)$.
\\ 
We also extend Stokes' formula to cut-off integrals of holomorphic families of symbol valued forms $\omega(z)$ obtained from a symbol valued  form  $\omega$ via a regularisation procedure (see Theorem \ref{thm:Stokesreg}):
$$\cutoffint d\,\left( \omega(z)\right)=0.$$
In the case of dimensional regularisation and when applied to forms built from Feynman type functions, this corresponds to a result already proven in \cite{E}.  
\section*{Acknowledgements} This paper was partially written during a
three months stay of the first author at the University Blaise Pascal in Clermont-Ferrand,  where
he was  invited on a  C.N.R.S.  position. The last author would like to
address her thanks to Edwin Langmann for pointing out to her reference
\cite{HH} as well as to Daniel Sternheimer for some interesting discussions
around this work. We realized at a late stage of the development of
  this work that our constructions are very close to prior results
of M. Lesch and M. Pflaum \cite{LP}, as a result of which we finally decided
not to publish this paper.  The last author thanks Mathias Lesch for drawing
her attention to his results with Markus Pflaum.\\
 However, our
emphasis is on symbols rather than on  operators so that  we are led to
considering cocycles that  mix the star and the ordinary function product on symbols  which do not
arise in their work.
\section{General scheme}
Take $X$ an open subset of $\R^m$, and
${\cal A} \subset C^\infty (X)$. Any associative (not necessarily
  commutative) product $\star$ on ${\cal A}$ induces a product $\wedge_\star$ on the  set  $\Omega {\cal A}$ of forms $\alpha$ on $X$ which are  of the type:
$$\alpha(x)=\sum_I \alpha_I (x) dx_I, \quad \alpha_I \in {\cal A}$$
as follows:
$$ (\alpha_I (x) dx_{i_1} \cdots dx_{i_p})\wedge_\star (\beta_J (x) dx_{j_1} \cdots dx_{j_q})=(\alpha_I \star \beta_J) (x) \, dx_{i_1} \cdots dx_{i_p} dx_{j_1} \cdots dx_{j_q}$$
which makes it a $\N$-graded algebra. If ${\cal A}$ is stable under partial derivations, then the exterior differential $d$ acts on  $\Omega{\cal A}$ increasing the degree by $1$:
$$d\,  (\alpha_I (x) dx_{i_1} \cdots dx_{i_p})= \sum_{j=1}^m\partial_j
\alpha_I (x)\, dx_jdx_{i_1} \cdots dx_{i_p}.$$
Clearly, equality $d^2=0$ comes from the odd parity of the $dx_i$
  which implies $dx_i dx_j= -dx_j dx_i$. \\
Let us furthermore assume that  partial derivations  $\partial_l$  on ${\cal A}$ obey the Leibniz rule,
i.e.  $$\partial_l
\left(a \star b \right)= \partial_l a\star b + a\star \partial_l b\quad
\forall \quad l=1, \cdots, m., \forall a, b\in {\cal A}.$$ Then $d$ is a graded derivation on $\Omega{\cal
  A}$; indeed, for any set of indices $I=\{i_1, \cdots, i_p\}$ and $J=\{i_{p+1}, \cdots,
i_{p+q}\}$  we have
\begin{eqnarray*} 
&{}& d\left(  (\alpha_I (x) dx_{i_1} \cdots dx_{i_p})\wedge_\star (\beta_J (x)
  dx_{i_{p+1}} \cdots dx_{i_{p+q}})\right)\\
&=&  d\left(  (\alpha_I (x) \star \beta_J (x) )\, dx_{i_1} \cdots dx_{i_p}
  dx_{i_{p+1}} \cdots dx_{i_{p+q}})\right)\\
&=&    \sum_{l=1}^{p+q} \partial_l
(\alpha_I (x)\star \beta_J (x))\, dx_l dx_{i_1} \cdots dx_{i_p}  \, dx_{i_{p+1}} \cdots dx_{i_{p+q}}\\
&=&    \sum_{l=1}^{p} \partial_l
\alpha_I (x)\star \beta_J (x)\, dx_l dx_{i_1} \cdots dx_{i_p}\, dx_{i_{p+1}} \cdots dx_{i_{p+q}}\\
&+&  (-1)^p \sum_{l=1}^{q}
\alpha_I (x)\star  \partial_l \beta_J (x)\,  dx_{i_1} \cdots dx_{i_p}\,dx_l
dx_{i_{p+1}} \cdots dx_{i_{p+q}}\\
&=&  d\left(  \alpha_I (x) dx_{i_1} \cdots dx_{i_p}\right)\wedge_\star (\beta_J (x)
  dx_{i_{p+1}} \cdots dx_{i_{p+q}})\\
&+& (-1)^p \left(\alpha_I (x) dx_{i_1} \cdots dx_{i_p}\right)\wedge_\star d\left(\beta_J (x)
  dx_{i_{p+1}} \cdots dx_{i_{p+q}}\right).
\\
\end{eqnarray*}

A linear map:
$\tau: {\cal A}\to \C$ induces a linear map $\bar \tau: \Omega {\cal A}\to \C$
defined by:
\begin{defn} Let $\alpha\in \Omega {\cal A}$. 
$$ \bar \tau\left( \alpha_I (x) dx_I\right)=\tau(\alpha_I), \quad {\rm
  if}\quad \vert I\vert =m; \quad \bar \tau\left( \alpha_I (x) dx_I\right)=0
\quad {\rm otherwise}.$$
We set: 
$$\chi_k^\tau(a_0, \cdots, a_k):= \bar \tau (a_0\star da_1\wedge_\star \cdots \wedge_\star
da_k).$$
\end{defn}
\begin{lem}
If $$\tau \left([a, b]_\star\right)=0\,\quad \forall a, b\in {\cal A}, \quad {\rm and} \quad \bar\tau \circ d=0$$
then $\left( \Omega {\cal A}, d, \bar \tau\right)$ defines  an
$m$-dimensional cycle with character  $(a_0, \cdots, a_m)\mapsto \chi_m^{\tau}(a_0, \cdots,
a_m)$ which yields 
 a cyclic Hochschild cocycle.
\end{lem}
{\bf Proof:} Since $\tau \left([a, b]_\star\right)=0$ we have that $$\tau \left(
  \alpha \wedge_\star \beta\right)= (-1)^{\vert\alpha\vert \cdot \vert \beta\vert}
\cdot \tau \left(\beta \wedge_\star \alpha\right),$$  which combined with $\bar \tau
\circ d=0$  provides an $m$-cycle. \endsquare
\begin{prop}
Let $\rho: {\cal A}\to \C$ be a linear map, let $\bar\rho: {\Omega\cal A}\to \C$
induced from $\rho$ as above, and let  $$\bar \tau:= \bar \rho
\circ d.$$ 
\begin{enumerate}
\item Then for any $a_1, \cdots, a_k\in {\cal A}$
$$B_0  \chi_k^\rho (a_1,
a_2,            \cdots ,             a_k)= \bar \tau ( a_1\star
da_2\wedge_\star \cdots \wedge_\star da_k).$$
\item If moreover there  is a trace $\tau$ on ${\cal A}$ (i.e. 
$\tau\left([a,b]_\star\right)=0\, \forall a, b\in {\cal A}$) such that
$\bar\tau$ coincides with the linear form on $\Omega{\cal A}$ associated with
$\tau$ as in definition 1 above, then 
$\left( \Omega {\cal A}, d, \bar \tau\right)$ defines  an
$m$-dimensional cycle with character  $$\chi_m^{\tau}(a_0, \cdots,
a_m):= \bar \tau (a_0\star da_1\wedge_\star \cdots \wedge_\star
da_m),$$ which yields 
 a cyclic Hochschild cocycle.
\end{enumerate}
\end{prop}
{\bf Proof}
\begin{enumerate}
\item Since $\bar \tau =\bar \rho \circ d$
\begin{eqnarray*}
B_0 \chi_k^\rho(a_1, \cdots, a_k)&=&  \chi_k^\rho(1, a_1, \cdots, a_k)\\
&=&  \bar \rho( da_1\wedge_\star \cdots \wedge_\star da_k)\\
&=& \bar \tau ( a_1\star da_2\wedge_\star \cdots \wedge_\star da_k).
\end{eqnarray*}
\item This follows from the above lemma since $\bar \tau \circ d= \bar \rho
  \circ d^2=0$.
\end{enumerate}
\endsquare
Let  ${\cal A}$ now be equipped with two (associative) products, the
  pointwise commutative one $\cdot$ and a non commutative one $\star$. Following  \cite{H} and \cite{CFS}
we set 
$$\theta(a, b):=  a\star b -a.b .$$
\begin{prop} Let $\tau: {\cal A}\to \C$ be  a trace with respect to
    the non-commutative product $\star$.
Then 
$$\phi^\tau_{2k}\left(a_0, a_1,\cdots, a_{2k}\right):= \tau \left(a_0\star \theta(a_1, a_2)\star  \cdots \star \theta(a_{2k-1}, a_{2k})\right)$$
defines a $b+B$-cocycle, namely
$$b \phi_{2k}+ \frac{1}{k+1}B \phi_{2k+2}=0.$$
\end{prop}
{\bf Proof:} 
The proof of \cite{H} and \cite{CFS} adapts to this general set up in a straightforward
manner. The assumption there that the star product be closed corresponds here
to the cyclicity of $\tau$.  One first shows that 
$\bar b \phi_{2k}=0$ with \begin{eqnarray*}
\bar b\chi(a_0,\cdots, a_{n})&=& \chi(a_0\star a_1, \cdots, a_j, a_{j+1}, \cdots, a_{n}) +\sum_{j=1}^{n-1} (-1)^j \chi(a_0, \cdots, a_j
\cdot a_{j+1}, \cdots, a_{n})\\
&+& (-1)^{n+1} \chi( a_0, \cdots,
a_{n-1}\star a_n).\end{eqnarray*}  
 The result then follows comparing $b\phi_{2k}$ and $\bar b \phi_{2k}$, which
 yields
\begin{eqnarray*}
b \phi_{2k}(a_0, \cdots, a_{2k+1})&=& -\tau\left( \theta(a_0, a_1)\star
  \theta(a_2, a_3)\star \cdots \right)+ \tau\left( \theta(a_{2k+1}, a_0)\star
  \theta(a_1, a_2)\star \cdots \right)\\
&=&-\frac{1}{k+1} B\phi_{2k+2}.
\end{eqnarray*}
\endsquare
\noindent In what follows we apply these constructions to the algebra of classical
symbols with compact support on an open subset of $\R^n$ letting $\tau$ be the Wodzicki residue on
 and $\rho$ the cut-off integral of symbols
\section{Classical symbols valued forms }
Let us first set some notations. \\
Let  $U$ be an open subset of $\R^n$. Let 
${\cal  S}^m(U)\subset \Ci(T^*U)$ denote  the set of
 scalar valued symbols on $U$ of order $m\in \R$, ${\cal S}(U):=\bigcup_{m\in
   \R}{\cal  S}^m(U)\subset C^\infty(T^*U)$ the algebra of all scalar valued
 symbols on $U$, ${\cal S}^{-\infty}(U):=\bigcap_{m\in \R}{\cal  S}^m(U)$ the
 algebra of scalar  smoothing symbols. We fix a norm on $\R^n$. Let
 $\chi$ be a smooth function on $T^*U$ such that $\chi(x,\xi)=0$ for $|\xi|\le
 1/2$ and $\chi(x,\xi)=1$ for $|\xi|\ge 1$.
\begin{defn}
$\sigma\in {\cal S}^m(U)$ is a classical symbol if  
for any positive integer $N$ we can write: 
\begin{equation}\label{eq:classicalsymb}
\sigma= \sum_{i=0}^{N} \chi \, \sigma_{m-i}+ \sigma_{(N)}.
\end{equation}
where $\sigma_{m-i}$ is positively homogeneous of order $m-i$ (i.e. 
 $$\sigma_{m-i}(x,t\xi)=t^{m-i} \sigma_{m-i}(x, \xi)$$
for any $t>0$ and any
 $(x, \xi) \in T_x^*U-\{0\}$), and where $\sigma_{(N)}$ is a symbol of order
 $m-N-1$.
We write for short 
$$\sigma\sim \sum_{i=0}^{\infty} \chi \, \sigma_{m-i}.$$
Let $CS^m(U)$ denote the class of scalar classical symbols of order $m$ and  
$CS(U)= \langle\bigcup_{m\in \C} CS^m(U)\rangle$ the algebra generated by
 scalar classical symbols of all orders. 
Similarly, let   $CS_{\rm com}^m(U)$ 
denote the subsets  of classical symbols of order $m$ with compact support in $U$ and $CS_{com}(U)= \langle\bigcup_{m\in \C} CS_{com}^m(U)\rangle$. \\
$CS^{\Z}(U):= \bigcup_{m\in \Z} CS^m(U)$ (resp. $CS_{com}^{\Z}(U):= \bigcup_{m\in \Z} CS_{com}^m(U)$) forms an alebra called the algebra of integer order symbols. We shall also consider its complement, namely the class  $CS^{\notin\, \Z}(U):=CS(U)- CS^{\Z}(U)$ (resp. $CS_{com}^{\notin\, \Z}(U):=CS_{com}(U)- CS_{com}^{\Z}(U)$) of non integer order symbols. 
 \end{defn}
Let us equip  $CS_{com}(U)$ with the left
product of symbols, also called the star product, which admits the following
asymptotic development:
$$(\sigma*\sigma')\sim\sum_{k\geq 0}(-i)^k\sum_{|\alpha|=k}
  {1\over \alpha!}\partial_\xi^\alpha\sigma.\partial_x^\alpha\sigma'.$$
(See for instance \cite{Sh} for details).
\\
\\
Symbol valued forms on $T^*U$ where $U$ is an open subset of $\R^n$ are defined as follows.
\begin{defn}Let $k$ be a non negative integer,
$m$ a complex number. We let
  \begin{eqnarray*}
\Omega^k\,CS^m(U)= \{\alpha \in \Omega^k(T^*U), 
\\
 \alpha &= &\sum_{I\subset \{1, \cdots, n\}, J\subset \{1, \cdots, n\}, \vert I
\vert + \vert J\vert =k} \alpha_{I, J} (x, \xi) \, dx_I 
\wedge d\xi_J\\
 {\rm with } &{}& \quad \alpha_{I, J} \in CS^{m-\vert J\vert}(U)\}
\end{eqnarray*}
denote the set of order $m$-classical symbol valued forms. 
  \end{defn}
The left product of symbols $\ast$ extends to symbol valued forms: given  
$$\alpha=  \sum_{I, J, \vert I
\vert + \vert J\vert =p} \alpha_{I, J} (x, \xi) \, dx_I \wedge d\xi_J\in
\Omega^pCS^m(U)$$ 
and 
$$\beta=  \sum_{K, L, \vert K
\vert + \vert L\vert =q} \alpha_{K, L} (x, \xi) \, dx_K \wedge d\xi_L\in
\Omega^qCS^n(U),$$
we set
\begin{equation}\label{eq:wedgeast}
\alpha\wedge_\ast \beta:=  \sum_{I, J, \vert I
\vert + \vert J\vert =k}\sum_{K, L, \vert K
\vert + \vert L\vert =q}  \alpha_{I, J} (x, \xi) \ast \alpha_{K, L} (x, \xi) \,   dx_I \wedge d\xi_J\wedge dx_K \wedge d\xi_L
\end{equation}
which lies in $\Omega^{p+q}CS^{m+n}(U)$. \\
Let $\Omega^k CS(U):= \langle \bigcup_{m\in\C} \Omega^k\,CS^m(U)\rangle $ (resp.  $\Omega^k CS_{com}(U):= \langle \bigcup_{m\in\C} \Omega^k\,CS_{com}^m(U)\rangle $) be the algebra generated by classical symbol (resp. with compact support)  valued $k$-forms of all orders. The sets $\Omega^k CS^{\Z}(U):=  \bigcup_{m\in\Z} \Omega^k\,CS^m(U) $, $\Omega^k CS_{com}^{\Z}(U):=  \bigcup_{m\in\Z} \Omega^k\,CS_{com}^m(U) $ form  algebras. We shall also consider the sets  
$\Omega^k CS^{\notin\, \Z}(U):=\bigcup_{m\notin\, \Z} \Omega^k\,CS^m(U) $ (resp. $\Omega^k CS_{com}^{\notin\, \Z}(U):=\bigcup_{m\notin\, \Z} \Omega^k\,CS_{com}^m(U) $). 
\begin{rk}
\begin{itemize}
\item With these conventions, $d\, \xi_j$ is of order $1$. Also, a $k$-form of order $0$ reads
$\alpha= \sum_{\vert I \vert + \vert J\vert =k} \alpha_{I, J} (x, \xi) \, dx_I\wedge d\xi_J$ with  $\alpha_{I, J}$ of order $-\vert J\vert$. 
\item The order of a zero degree symbol valued form $\sigma\in \Omega^0\,
  CS^m(U)$ coincides with the order of the corresponding  classical symbol
  $\sigma$.
\item More generally, any  zero order symbol valued $k$-form on $U$ is of the type 
$$\alpha=\sum_{\vert I\vert +\vert J\vert=k}\alpha_{I, J} d\, x_I\wedge d\,
\xi_J  $$ with $\alpha_{I, J}$ of order $-\vert J\vert$. In particular, given
any $\sigma\in CS(U)$, the top form $\sigma_{-n}(x,\xi)\,dx_1\wedge\cdots\wedge dx_n\wedge d\xi_1\wedge\cdots \wedge d\xi_n$
 provides an example of positively homogeneous zero order symbol valued
 $n$-form. 
\end{itemize}
\end{rk}
\begin{lem}
A classical symbol valued form $\alpha\in \Omega^k\, CS^m(U)$  of order $m$ has an asymptotic expansion of the following form. For any non negative integer $N$, there is a symbol valued form  $\alpha_{(N)}$ of order $m-N-1$ such that
$$\alpha= \sum_{i=0}^N \alpha_{m-i}+ \alpha_{(N)}$$
with  $\alpha_{m-i}:=  \sum_{\vert I\vert+\vert J\vert=k}\alpha_{I, J,m-\vert J\vert - i}\, dx_I\wedge d\xi_J$ is positively homogeneous of order $m-i$, with $\alpha_{I, J,m-\vert J\vert - i}$ positively homogeneous of order $m-\vert J\vert -i$. \\
Furthermore, the exterior differentiation $d$ sends  $\Omega^k\, CS^m(U)$ to $ \Omega^{k+1}\, CS^m(U)$ and  for any  integer $j\leq m$, we have 
 $$\left(d\alpha\right)_j=d\, \alpha_j.$$ 
\end{lem}
{\bf Proof:}
The first part of the statement follows trivially from the description of
$\alpha$ combined with the properties of ordinary classical symbols. As for the second   part of the statement we write
\begin{eqnarray} \label{eq:dalpha}
 d\, \alpha&=& d\, \sum_{\vert I\vert+\vert J\vert=k}\alpha_{I, J}\, dx_I\wedge d\xi_J\nonumber\\
&=&  \sum_{l=1}^{n}\sum_{\vert I\vert+\vert J\vert=k}
\frac{\partial}{\partial x_l}\alpha_{I, J}\,dx_l\wedge dx_I\wedge d\xi_J\nonumber\\
&+&   \sum_{l=1}^{n}\sum_{\vert I\vert+\vert J\vert=k}\frac{\partial}{\partial \xi_m}\alpha_{I, J}\,d\xi_m\wedge dx_I\wedge d\xi_J
\end{eqnarray}
which lies in $\Omega^k CS(U)$ since the order of $\frac{\partial}{\partial
  x_m}\alpha_{I, J} \, d\xi_m$ coincides with that of $\alpha_{I, J}$. The
computation above also shows that if $\alpha$ is positively homogeneous of
order $m$, so is $d\,\alpha$, which ends the proof of the lemma.
 \endsquare
\begin{rk}
In particular, for $\alpha\in \Omega\, CS(U)$ we have:
$$\left(d\alpha\right)_{0}= d\, \alpha_{0}.$$
\end{rk}
\section{The Wodzicki residue character  on  classical symbols }
\subsection{The Wodzicki residue extended to classical symbol valued valued 
forms}
Let us first briefly  recall the notion of Wodzicki residue on classical
symbols \cite{W}, \cite{K}.
\begin{defn} Let  $U$ be  an open subset in $\R^n$ and $x$ a point in $U$. The
  (local) Wodzicki residue density of a classical
  symbol $\sigma \in  CS(U)$ at point $x$  is given by $$ {\rm
res}_x(\sigma)=\int_{\vert \xi\vert=1}\sigma_{-n} (x,\xi) \, d_S\xi, $$
where $d_S\xi=\sum_{i=1}^n (-1)^{i+1} \xi_i\,d\xi _1\wedge \cdots \wedge d\hat \xi_i\wedge \cdots \wedge d\xi_n $  and $\vert \xi \vert= (\sum_{i=1}^n \xi_i^2)^{1/2}$ is the canonical norm in $\R^n$. \\
For any $\sigma\in CS(U)$ with compact support the Wodzicki residue of
$\sigma$ is then
defined as:
$${\rm res}(\sigma):=\int_U {\rm res}_x(\sigma) \, dx.$$
\end{defn}
\begin{rk} For any $t>0$ we have $d_S(t\xi)= t^n d_S\xi$ and $\sigma_{-n}
  (x,t\xi)= t^{-n} \sigma_{-n} (x,\xi)$ so that the form $\sigma_{-n} (x,\xi)
  \, dx\wedge d_S\xi$ is positively homogeneous of degree $0$. 
\end{rk}
The Wodzicki residue extends  from $CS(U)$ to
$\Omega\,CS(U)$ in a straight forward manner.
\begin{defn} For any $\alpha= \sum_{I, J} \alpha_{IJ} \, dx_I\wedge d\xi_J\in
\Omega CS_{\rm com}(U)$, for any $x\in U$ we set
$${\rm res}_x\left(\sum_{I, J} \alpha_{IJ} \, dx_I\wedge d\xi_J\right)=\sum_I {\rm res}_x
  (\alpha_{IJ})\, dx_I=\sum_I \int_{\vert \xi\vert=1} \left(\alpha_{IJ}\right)_{-n} (x, \xi)
  \, d_S \xi, \quad {\rm if} \quad \vert J\vert =n$$
and ${\rm res}_x\left(\sum_{I, J} \alpha_{IJ} \, dx_I\wedge d\xi_J\right)=0$
whenever   $\vert J\vert\neq n$.\\ Similarly, we set: 
$${\rm res}\left(\sum_{I, J} \alpha_{IJ} \, dx_I\wedge d\xi_J\right)= {\rm res}
  (\alpha_{IJ})=\int_{\vert \xi\vert=1} \left(\alpha_{IJ}\right)_{-n} (x, \xi)
  \, dx\wedge d_S \xi, \quad {\rm if} \quad \vert I\vert =\vert J\vert=n$$
and ${\rm res}\left(\sum_{I, J} \alpha_{IJ} \, dx_I\wedge d\xi_J\right)=0$
whenever $\vert I\vert\neq n$ or  $\vert J\vert\neq n$. 
 \end{defn}
It is useful to give an alternative more intrinsic formulation of this
extended Wodzicki residue. The form $d_S\xi$ on $T^*_xU$ can be seen as the interior product $i_X(\Omega_x)$ of the volume form $ \Omega_x:= d\xi_1\wedge\cdots\wedge d\xi_n$ on $T^*_x U$ with the Liouville (or radial) field $$X(x, \xi)= \sum_{i=1}^n \xi_i \frac{\partial}{\partial \xi_i}.$$  This Liouville  field can also be seen as the generator 
$$ X(x, \xi):= \frac{d}{dt}_{\vert_{t=0}}f_t(x, \xi) $$ of the one parameter semigroup of transformations of  $T^*U$:
\begin{eqnarray*}
 \R\times T^*U &\to & T^*U\\
\left(t, (x, \xi)\right) &\mapsto & f_t(x, \xi):= (x, e^t\, \xi).\\
\end{eqnarray*} 
Let $\rho: T^*U-\{0\}\to S^*U$ denote the radial projection
 $\rho(x, \xi)= (x, \frac{\xi}{\vert \xi\vert})$, and let $j:S^*U\to
 T^*U-\{0\}$ denote the canonical fibre bundle injection. Clearly
 $\rho\circ j=Id$. We have the following lemma.
\begin{lem}
A form $\alpha$ on $T^*U-\{0\}$ is positively homogeneous of order zero if and only if it satisfies one of the two equivalent conditions:
\begin{enumerate}
\item the form  can be written
\begin{equation}\label{eq:0form}
\alpha = \rho^*\beta+{dr\over r}\wedge \rho^*\gamma
\end{equation}
with $\beta,\gamma\in\Omega(S^*U)$, and more precisely:
$$\beta=j^*\alpha,\hskip 20mm \gamma=j^*(\iota_X\alpha).$$
\item 
${\cal L}_X(\alpha)= 0$
where ${\cal L}_X$ is the Lie derivative in direction $X$.
\end{enumerate}
\end{lem}
{\bf Proof:} The second condition is equivalent to 
$\alpha(x, e^t \, \xi)= \alpha(x, \xi) \quad \forall t>0$ and hence to positive homogeneity of order zero 
since:
$${\cal L}_X\alpha=\frac{d}{dt}_{\vert_{t=0}}  f_t^* \alpha=
\frac{d}{dt}_{\vert_{t=0}} \alpha(x, e^t \, \xi).$$
For any
$\beta\in\Omega(S^*U)$ the differential form
$\rho^*\beta$ is invariant by dilations, hence positively homogeneous of
order zero. The first condition then clearly
implies that $\alpha$ is positively homogeneous of order zero, as ${dr\over
  r}$ obviously is, hence $(1)\Rightarrow (2)$. Suppose now that (2) is
verified, and seek for $\beta$ and $\gamma$ such that \ref{eq:0form}
holds. As $j^*({dr\over r})=0$ and $\rho\circ j=Id$ we clearly
have:
$$j^*\alpha=j^*(\rho^*\beta+{dr\over r}\wedge \gamma)=\beta.$$
Now $\iota_X\alpha=\iota_X{dr\over r}\wedge\rho^*\gamma =\rho^*\gamma$, hence
$\gamma=j^*\rho^*\gamma=j^*(\iota_X\alpha)$. We have then proved the
uniqueness of $\beta$ and $\gamma$. To prove the existence, notice that the
difference:
$$\delta=\alpha-(\rho^*j^*\alpha+{dr\over r}\wedge\rho^*j^*\iota_X\alpha)$$
verifies $j^*\delta=\iota_X\delta=0$, hence it easily follows that $\delta=0$.
So $(2)\Rightarrow (1)$. \endsquare 
\vskip 3mm\noindent
\begin{ex}Given any $\sigma\in CS(U)$, the top form 
 \begin{eqnarray*}
\alpha_\sigma(x, r\cdot \omega)&:=&\sigma_{-n}(x,r\cdot
\omega)\,dx_1\wedge\cdots\wedge dx_n\wedge \frac{dr}{r}\wedge d_S\omega\\
&= &\sigma_{-n}(x,\xi)\,dx_1\wedge\cdots\wedge dx_n\wedge d\xi_1\wedge\cdots \wedge d\xi_n
\end{eqnarray*}
is a  positively homogeneous zero order symbol valued
 $n$-form and we have:
$$\iota_X\alpha_\sigma=\sigma_{-n}dx_1\wedge\cdots\wedge dx_n\wedge d_S\xi.$$
\end{ex}
The following elementary result provides a more intrinsic formulation of the
Wodzicki residue extended to forms. 
\begin{prop}
Let $U$ be an open subset of $\R^n$. Denote by $j$ (resp. $j_x$ for any $x\in U$) the
injection of $S^*U$ (resp. $S_x^*U$) inside the cotangent bundle $T^*U$
 (resp. inside $T_x^*U$).
Given $\alpha\in \Omega\, CS(U)$,  for any $x\in U$:
$${\rm res}_x (\alpha):=
\int_{S_x^*U} j_x^*(\iota_\Lambda\iota_X\alpha_0)$$
where $\Lambda$ stands for the volume element $n!{\partial\over\partial x_1}
\wedge\cdots\wedge{\partial\over\partial x_n}$, and:
$${\rm res} (\alpha):=\int_U{\rm res}_x(\alpha)\, dx_1\cdots dx_n=
\int_{S^*U} j^*(\iota_X\alpha_0).$$
\end{prop}
\subsection{Stokes' formula for the Wodzicki residue} 
 \begin{thm}\label{thm:Stokesres}
For  any $\beta\in  \Omega\, CS(U) $ with compact support we have
$${\res}\left( d\beta\right)=0.$$
  \end{thm}
{\bf Proof: } Using Cartan's formula, this follows from Stokes' property for ordinary integrals, since $\left(d\beta\right)_0= d\,
\beta_0$ implies ${\cal L}_Xd\beta_0=0$, hence: 
$${\rm res}\,(d\beta)= \int_{S^* U}j^*(\iota_X d\beta_0)= -\int_{S^*
  U}j^*(d\iota_X\beta_0)=-\int_{S^* U}d\,\bigl(j^*(\iota_X\beta_0)\bigr)=0$$
since  $S^*U$ is boundaryless.  \endsquare
We recover this way a known integration by parts formula for the Wodzicki residue which underlies the traciality property of the Wodzicki residue on classical pseudodifferential  operators.  
\begin{cor}
 For any  $\sigma\in  CS(U) $ with compact support, 
$${\rm res}\left(\frac{\partial }{\partial \, \xi_i}\sigma\,  \sigma^\prime \right)=-{\rm res}\left(\sigma\, \frac{\partial }{\partial \, \xi_i} \sigma^\prime \right) \quad \forall i\in \{1, \cdots, n\}$$
and 
$${\rm res}\left( \frac{\partial }{\partial \, x_i}\sigma\,  \sigma^\prime \right)=-{\rm res}\left( \sigma\, \frac{\partial }{\partial \, x_i} \sigma^\prime \right)\quad \forall i\in \{1, \cdots, n\}.$$
\end{cor}
{\bf Proof:} Let  $\tau\in CS(U)$ with compact support.  Applying Theorem \ref{thm:Stokesres} to $\beta_\tau^i:= \tau_{-n+1}(x, \xi) \,  d\xi_1\wedge \cdots \wedge \hat{d\xi_i}\wedge \cdots \wedge d\xi_n$ we get 
\begin{eqnarray*}
{\rm res}\left( \frac{\partial }{\partial \, \xi_i}\tau(x, \xi) \right)&=& {\rm res}\left( \left(\frac{\partial }{\partial \, \xi_i}\tau(x, \xi)\right)_{-n}\,dx_1\wedge\cdots\wedge dx_n\wedge d\xi_1\wedge\cdots d\xi_n \right)\\
&=& {\rm res}\left( \frac{\partial }{\partial \, \xi_i}\tau_{-n+1}(x, \xi)\,dx_1\wedge\cdots\wedge dx_n\wedge d\xi_1\wedge\cdots d\xi_n \right)\\
&=&(-1)^{i-1}{\rm res}\left(d\left(\tau_{-n+1}(x, \xi) \,  d\xi_1\wedge \cdots \wedge \hat{d\xi_i}\wedge \cdots \wedge d\xi_n\right)\right)\\
&=&0.
\end{eqnarray*}
Applying this to $\tau := \sigma\, \sigma^\prime$ yields the first part of the corollary. 
A  similar proof  replacing $\frac{\partial }{\partial \, \xi_i}$ by $\frac{\partial }{\partial \, x_i}$ using Stokes' formula applied to  $\beta_\tau^i:= \tau_{-n}(x, \xi) \,  d\xi_1\wedge\cdots\wedge d\hat x_i
\wedge  d\xi_1\wedge \cdots \wedge d\xi_n$ gives the second equality of the corollary. 
\endsquare
\begin{cor} The Wodzicki residue defines a trace on the subalgebra $CS_{\rm com}(U)\in CS(U) $ of symbols with compact support in $x$
$${\rm res}([\sigma, \sigma^\prime ]_\ast)=0\quad\forall \sigma, \sigma^\prime \in \Cl_{\rm comp}(U,\C)$$
where we have set
$[\sigma, \sigma^\prime ]_\ast:= \sigma\ast \sigma^\prime -\sigma\ast \sigma^\prime .$
\end{cor}
{\bf Proof:}  We use the asymptotic development of the left product of
symbols. There exists a positive integer $N$ such that~:
$${\rm res}(\sigma\ast\sigma')=\sum_{k\le N}(-i)^k\sum_{|\alpha|=k} 
{1\over\alpha!}{\rm res}(\partial_\xi^\alpha\sigma.\partial_x^\alpha\sigma').$$
Indeed, the remainder term will be of order $<-n$ for sufficiently big $N$,
and then will have vanishing residue. By the above lemma, we have for $\sigma,\sigma^\prime\in CS(U)$ with compact support in $U$
\begin{eqnarray*}
{\rm res}(\sigma\ast\sigma^\prime)&=& \sum_{\vert
\alpha\vert\le N}i^{|\alpha|}\frac{1}{\alpha!}{\rm res}\left( \partial_\xi^\alpha \sigma\cdot
\partial_x^\alpha \sigma^\prime\right)\\
&=& \sum_{\vert
\alpha\vert\le N}(-i)^{|\alpha|}\frac{1}{ \alpha!}{\rm res}\left(\partial_x^\alpha \partial_\xi^\alpha \sigma\cdot
 \sigma^\prime\right)\\
&=& \sum_{\vert
\alpha\vert\le N}i^{|\alpha|}\frac{1}{ \alpha!}{\rm res}\left(\partial_x^\alpha  \sigma\cdot
\partial_\xi^\alpha \sigma^\prime\right)\\
&=& {\rm res}(\sigma^\prime \ast\sigma).
\end{eqnarray*} 
\endsquare
\subsection{A Wodzicki residue cycle on zero order classical symbols}
The exterior  differential:
$$d: \Omega^k\, CS(U)\to \Omega^{k+1}\, CS(U)$$
obeys the usual ``Leibniz rule'': 
$$d\left(\alpha\wedge_\ast \beta\right)= d\alpha\wedge_{\ast}\beta+ (-1)^k\alpha\wedge_{\ast}d\beta \quad\forall 
\alpha\in\Omega^kCS(U), \beta\in \Omega^*CS(U)$$
as can easily be checked from (\ref{eq:wedgeast}) and (\ref{eq:dalpha}) so that $\left(\Omega\, CS_{\rm com}(U), d\right)$ is a graded differential algebra  with  $CS_{\rm com}(U)$  equipped with the left product of  symbols. 
\begin{thm}\label{thm:residuecycle}
Let $CS_{\rm com}(U)$ be equipped with the left product of  symbols.  The triple $\left(\Omega\, CS_{\rm com}(U), d, {\rm res}\right)$ yields an $2n$-cycle which we refer to as the Wodzicki residue cycle. 
\end{thm}
{\bf Proof:} As previously observed, the Wodzicki residue vanishes on $\Omega^kCS_{com}(U)$ for $k<2n$. It is closed by the Stokes' formula since ${\rm res}(d\beta)=0$ for any  $\beta\in \Omega^{2n-1}CS_{com}(U)$.\\
The fact that the ordinary Wodzicki residue defines a trace on $CS_{\rm
  com}(U)$ immediately implies:
$${\rm res}(\alpha\wedge_\ast \beta)= (-1)^{\vert \alpha\vert \cdot \vert \beta\vert} {\rm res} (\beta\wedge_\ast \alpha)$$
so that $\left(\Omega\, CS_{\rm com}(U), d, {\rm res}\right)$ defines a cycle.\endsquare
We call {\bf residue character} the  associated $2n$-character (see Appendix A). 
\begin{defn}
 Let the {\rm residue} $k$-cochain denote the $k+1$-linear form on $CS_{\rm com}(U)$ 
\begin{eqnarray*}
\chi_{k}^{\rm res}(\sigma_0, \cdots, \sigma_{k})&=& {\rm res}\left( \sigma_0\ast d\, \sigma_1\wedge_\ast \cdots \wedge_\ast d\,\sigma_{k}\right)
\end{eqnarray*}
for all $ \sigma_0, \cdots, \sigma_{k} \in CS_{\rm com}(U)$.
\end{defn}
Residue $k$-cochains vanish for $k<2n$ and 
the residue character is the $2n$-residue cochain $\chi_{2n}^{\rm res}$. It
satisfies the following properties (with the notations of Appendix A):
\begin{itemize}
\item   $B_0\chi_{2n}^{\rm res}=0$ and $B\chi_{2n}^{\rm res}=0$, 
\item  $b_\ast\chi_{2n}^{\rm res}=0$ where $b_\ast$ is the Hochschild coboundary operator associated with the left product on symbols.
\end{itemize} 
Restricting to zero order symbols we get:
\begin{thm}
For any symbols  $  \sigma_0, \cdots, \sigma_{2n} \in CS^0_{\rm
    com}(U),$
\begin{eqnarray*} 
\chi_{2n}^{res}(\sigma_0, \cdots,
\sigma_{2n})&=&\rm{res}(\sigma_0\, d\sigma_1\wedge\cdots\wedge d\sigma_{2n})\\
&=&\int_{S^*U} j^*\iota_X(\sigma_0^L\, d\sigma_1^L\wedge\cdots\wedge
d\sigma_{2n}^L)\\
&=&\frac{(-1)^n}{n!}\, A\left[ {\rm res}\left(\sigma_0\, \theta \left(\sigma_1,\sigma_2\right)\cdots\, \theta\left(\sigma_{2n-1},\sigma_{2n}\right)\right)\right].\\
 \end{eqnarray*}
Here $\sigma_i^L$ stands for the leading symbol of $\sigma_i$ and where we
have set $\theta(\sigma_i, \sigma_j)= \sigma_i \star \sigma_j - \sigma_i\cdot
\sigma_j$  as in section 1.  $A$ denotes the antisymmetrisation over all
but the first variable. 
\end{thm}
{\bf Proof:} The difference $\sigma_0*d\sigma_1\wedge_\ast\cdots \wedge_\ast d\sigma_{2n}
-\sigma_0\,d\sigma_1\wedge\cdots\wedge d\sigma_{2n}$ has clearly vanishing residue as top
form of order $\le -1$, hence the first equality. The second equality then follows since the top order term
$(\sigma_0d\sigma_1\wedge\cdots\wedge d\sigma_{2n})_0$ is precisely
$\sigma_0^Ld\sigma_1^L\wedge\cdots\wedge d\sigma_{2n}^L$. \\
As for the last equality, we have $\theta(\sigma_i, \sigma_j)\sim \sum_{\vert
  \alpha\vert \neq 0 }\frac{(-i)^{\vert  \alpha\vert}}{\alpha!} \partial^\alpha \sigma_i \partial_x^\alpha
  \sigma_j$
so that 
\begin{eqnarray*}
&{}&A\left[ {\rm res}\left(\sigma_0\, \theta \left(\sigma_1,\sigma_2\right)\cdots\,
  \theta\left(\sigma_{2n-1},\sigma_{2n}\right)\right)\right]\\
&=&A\left[ \sum_{\vert
  \alpha_1\vert \neq 0,\cdots,\vert
  \alpha_n\vert\neq 0 } \frac{(-i)^{\vert
    \alpha\vert}}{\alpha_1!\cdots\alpha_n!}\, {\rm
    res}\left(\sigma_0\, \partial_\xi^{\alpha_1}\sigma_1\, \partial_x^{\alpha_1}\sigma_2\cdots\,
  \partial_\xi^{\alpha_n}\sigma_{2n-1}\,
  \partial_x^{\alpha_n}\sigma_{2n}\right)\right]\\
&=&    A\left[\sum_{i_1,\ldots ,i_n\in\{1,\ldots,n\}} (-i)^n{\rm
    res}\left(\sigma_0\, \partial_{\xi_{i_1}}\sigma_1\, \partial_{x_{i_1}}\sigma_2\cdots\,
  \partial_{\xi_{i_n}}\sigma_{2n-1}\,
  \partial_{x_{i_n}}\sigma_{2n}\right)\right]\\
&=&i^nn!  \, {\rm
res}\left(\sigma_0\,d\sigma_1\wedge d \sigma_2\wedge \cdots\wedge
  d\sigma_{2n-1}\wedge
  d\sigma_{2n}\right).\\
\end{eqnarray*}
\endsquare

\section{Cut-off integrals of symbol valued forms and cosphere cochain}
\subsection{Cut-off integrals extended to classical symbols valued forms} 
Defining cut-off integrals amounts to extracting   finite parts from otherwise divergent
integrals, a procedure which we recall here (without proofs) in the case of ordinary classical symbols \cite{H}, \cite{G}, \cite{W}, \cite{KV}. 
 \begin{prop} \label{prop:cutoffint} Let $U$ be an open subset of $\R^n$ and let $x\in U$. Given $\sigma\sim \sum_{i=0}^\infty \chi\,
\sigma_{m-i} \in CS^m(U)$, the expression $\int_{B_x^*(0,R)} \sigma(x,\xi)\, d\xi$
has an asymptotic expansion  $$\int_{B_x^*(0, R)} \sigma(x,\xi) \, d\xi=c(x)+
\sum_{i=0, \,m-i+n\neq 0 }^\infty  a_{i}(x) \,\frac{R^{m-i+n}}{m-i+n}+
b(x)\log R $$ where $c(x), a_i(x), b(x)\in\C$.
The finite part   called the cut-off integral of $\sigma(x, \cdot)$ which is  given by the constant  $c(x)$ reads:
\begin{eqnarray}\label{eq:finitepart} \cutoffint_{T_x^*U} \sigma(x,\xi) \,
d\xi&:=& {\rm fp}_{R\to\infty} \int_{B_x^*(0, R)} \sigma(x,\xi) \, d\xi\\
&=& \int_{B_x^*(0,1)} \sigma(x,\xi)\,d\xi\nonumber\\ &\hskip 3mm +&
\int_{T_x^*U-B_x^*(0, 1)} \sigma_{(N)} (x, \xi) \, d\xi\nonumber \\ &\hskip 3mm -
&\sum_{i=0, \, m-i+n\neq 0 }^{K_N}\, \frac{1}{m-i+n} \,
\int_{\vert \xi \vert=1}\sigma_{m-i}(x,\xi) \, d_S\xi.
\end{eqnarray}  \\ 
If $\sigma\in CS_{com}(\rm U)$ we set
$$ \cutoffint_{T^*U} \sigma(x,\xi) := \int_U dx\, \cutoffint_{T_x^*U} \sigma(x,\xi) .$$
 The constant  
$ b(x)$ coincides with the local Wodzicki residue density
${\rm res}_x(\sigma)$. When it vanishes, the  finite part ${\rm fp}_{R\to\infty} \int_{B_x^*(0, R)} \sigma(x,\xi) \, d\xi$ is independent of the rescaling $R\mapsto \lambda R$. Specifically, this holds for non integer order symbols. 
\end{prop}
\begin{rk}
This cut-off integral extends the ordinary integral in the following sense; if $\sigma$ has order smaller than $-n$ then $\int_{B_x^*(0, R)} \sigma(x,\xi)\, d\xi$ converges when $R\to \infty$ and $\cutoffint_{T_x^*U} \sigma(x,\xi)\, dx=\int_{T_x^*U} \sigma(x,\xi) \, d\xi.$
\end{rk}
\begin{defn}The cut-off integral on $T_x^* U$  of a form $\alpha=\sum_{I,
    J}\alpha_{I, J} dx_I \wedge d\xi_J\in
  \Omega\, CS(U) $ with compact support in $x$ is defined by: 
$$\cutoffint_{T^*U}
\alpha:= \cutoffint_{T^*U}\alpha_{I, J}(x, \xi) \, dx_1\wedge\cdots\wedge dx_nd\xi_1\wedge d\xi_n\quad {\rm if} \quad \vert I \vert= \vert J \vert=n$$
and which vanishes otherwise.
\end{defn}
As in the case of ordinary integrals, we
recover the cut-off integral on symbol valued functions $\sigma\in CS^{\notin
\, \Z}(U) $ via the integral on forms  by integrating the top form $\sigma(x,
\xi) \, dx\wedge d\xi$ setting:
$$\cutoffint_{T^*U} \sigma(x, \xi):= \cutoffint_{T_x^*U} \sigma(x, \xi) \,
 dx\wedge d\xi$$ where the right hand side is now seen as a cut-off integral on a symbol valued form.\\ \\
Similarly to   ordinary integrals, cut-off integrals on forms satisfy Stokes'
property (compare with  Lemma 5.5 in \cite{LP}).
\begin{thm}\label{thm:Stokescutoff}
Let $U$ be an open subset of $\R^n$ and let $\beta\in \Omega^{2n-1}\, CS_{com}(U) $ be a symbol valued form. Then 
$$\cutoffint_{T^*U} d\beta=\sum_{I, J}\,  \int_{S^*(0, 1)}  \beta_{I,J,\, -n+1}(x, \xi)\,  dx_I \wedge    d\xi_J$$
so that Stokes' formula:
$$\cutoffint_{T^*U} d\beta=0$$ holds whenever $\beta \in\Omega^{2n-1}\, CS_{com}^{\notin \, \Z}(U) $. \\
Here $\beta(x,\xi)=   \sum_{I, J\subset \{1, \cdots, n\}, \vert I\vert +\vert J\vert=2n-1} \beta_{I, J}(x, \xi)\, dx_I \wedge d\xi_J $ with $\beta_{I, J}\in CS(U)$.
\end{thm}
{\bf Proof:}
The  $2n-1$ form reads  $\beta(x,\xi)=   \sum_{I, J\subset \{1, \cdots, n\}, \vert I\vert +\vert J\vert=2n-1} \beta_{I, J}(x, \xi)\, dx_I \wedge d\xi_J $ with $\beta_{I, J}\in CS^{\notin\, \Z}(U)$ so that, letting $B^*(0, R)$, resp. $S^*(0, R)$ be respectively the ball in the cotangent bundle of radius $R$ centered at the origin, and the sphere in the cotangent bundle of radius $R$ centered at the origin, we have
\begin{eqnarray*}
\cutoffint_{T^*U} \, d\beta&=&   \sum_{I, J}\cutoffint_{T^*U}  d\,\left(
\beta_{I, J}(x, \xi)\,  dx_I\wedge d\xi_J\right)\\
&=&\sum_{I, J}{\rm fp}_{R\to \infty}  \int_{B^*(0, R)} d\, \left(\beta_{I, J}(x, \xi)\,  d x_{I} \wedge    d\xi_{J}\right)\\
&=&\sum_{I, J}{\rm fp}_{R\to \infty} \int_{S^*(0, R)} \beta_{I, J}(x, \xi)\,  dx_I \wedge    d\xi_{J}\\
&{}& {\rm using} \quad {\rm Stokes'} \quad {\rm property} \quad {\rm for} \quad
 {\rm ordinary} \quad {\rm integrals}\\
&=& \sum_{I, J}\sum_{j=0}^N\, {\rm fp}_{R\to \infty} \,
\int_{S^*(0, R)} \chi(\xi)\, \beta_{I, J, \,m_{I,J}-j-\vert J\vert}(x, \xi)\,  d x_I\wedge    d\xi_J\\
&\hskip 12mm +& \lim_{R\to \infty}\int_{S^*(0,R)}\beta_{I,J, \, (N)} \\
&{}& \hskip 12mm ({\rm where} \quad \beta_{I,J}= \sum_{j=0}^{N} \chi\,\beta_{I,J,\,m_{I,J}-j-\vert J\vert}+ \beta_{I,J, \, (N)})\\ 
&= &  \sum_{I, J}\sum_{j=0}^{N}\, {\rm fp}_{R\to \infty} \,
\int_{S^*(0, R)} \beta_{I,J, \,m_{I,J}- j-\vert J\vert}(x, \xi)\,  dx_I  \wedge d\xi_J \\
&{}& \hskip -9mm{\rm since } \ \lim_{\vert \xi \vert \to \infty}|\xi|^{n-1}\beta_{I,J, \,
  (N)}(x, \xi)=0\ {\rm and}\ \chi=1\  {\rm outside}  \ B^*(0, 1)\\
&= &  \sum_{I, J}\sum_{j}^{N}\, {\rm fp}_{R\to \infty} \,
 R^{m_{I,J}-j-\vert J\vert +n-1}\,  \int_{S^*(0, 1)}  \beta_{I,J,\, m_{I,J}-j-\vert J\vert}(x, \xi)\,  dx_I \wedge    d\xi_J \\
&=&0\quad {\rm whenever} \quad  m_{I,J}-j-\vert J\vert+n-1\neq 0,\\
&=& \sum_{I, J}\,  \int_{S^*(0, 1)}  \beta_{I,J,\, -n+1}(x, \xi)\,  dx_I \wedge    d\xi_J\quad {\rm otherwise} \\
\end{eqnarray*}
where $m_{I,J}\notin \Z$ is the order of $\beta_{I,J}$.\endsquare
 As a consequence, cut-off integrals on non-integer order  symbols satisfy an integration by parts formula:
\begin{cor}  
For  any  $\sigma\in  CS_{com}(U) $ then 
$$\cutoffint_{T^*U} \frac{\partial }{\partial \, \xi_i}\sigma(x, \xi) \, d\xi\, dx=
(-1)^{i-1}\int_{S^*(0, 1)} \sigma_{-n+1}(x, \xi) \,  d\xi_1\wedge \cdots \wedge \hat{d\xi_i}\wedge \cdots \wedge d\xi_n\wedge dx_1\wedge\cdots \wedge dx_n.$$
 In particular, if  $\sigma\in  CS_{com}^{\notin \, \Z}(U) $ then 
$$\cutoffint_{T^*U} \frac{\partial }{\partial \, \xi_i}\sigma(x, \xi) \, d\xi\, dx=0\quad \forall i\in \{1, \cdots, n\}.$$
\end{cor}
{\bf Proof:}
  Applying 
 Stokes' formula to $\beta:= \sigma(x, \xi) \,  d\xi_1\wedge \cdots \wedge \hat{d\xi_i}\wedge \cdots \wedge d\xi_n \wedge dx_1\wedge\cdots\wedge dx_n$ we have:
\begin{eqnarray*}
\cutoffint_{T^*U} \frac{\partial }{\partial \, \xi_i}\sigma(x, \xi) \,  d\xi&=& (-1)^{i-1}\cutoffint_{T^*U} d\left(\sigma(x, \xi) \,  d\xi_1\wedge \cdots \wedge \hat{d\xi_i}\wedge \cdots \wedge d\xi_n\wedge dx_1\wedge\cdots \wedge dx_n\right)\\
&=& (-1)^{i-1}\int_{S^*(0, 1)} \sigma_{-n+1}(x, \xi) \,  d\xi_1\wedge \cdots \wedge \hat{d\xi_i}\wedge \cdots \wedge d\xi_n\wedge dx_1\wedge\cdots \wedge dx_n .
\end{eqnarray*}
This last term vanishes whenever $\sigma$ has non integer order. \endsquare
The integration by parts formula  yields translation invariance of cut-off integrals on non integer order symbols.
\begin{cor}For any $\sigma\in CS_{com}(U)$
$$\cutoffint_{T^*U} \sigma(x,\xi+\eta)\,dx\,  d\xi=\cutoffint_{T^*U}  \sigma(x,\xi)\,dx\,  d\xi\quad \forall \eta\in \Ci(U, T^* U).$$
If $\sigma\in CS_{com}^{\notin\, \Z}(U)$ then  
$$\cutoffint_{T^*U} \sigma(x,\xi+\eta)\,dx\,  d\xi=\cutoffint_{T^*U}  \sigma(x,\xi)\,dx\,  d\xi\quad \forall \eta\in \Ci(U, T^* U).$$
 \end{cor}
{\bf Proof:} A Taylor expansion $\eta\mapsto \sigma(\xi+ \eta)$ in $\eta$ at $0$ yields, for any $x\in U$,  the existence of some $\theta\in ]0, 1[$ such that: 
$$
\cutoffint_{T_x^*U} \sigma(x,\xi+\eta)\, d\xi= \sum_{\vert \alpha\vert \leq K}\cutoffint_{T_x^*U} d\xi\, \frac{ D_\xi^\alpha\sigma(x,\xi)}{\alpha!} \eta^\alpha+ \sum_{\vert \alpha\vert = K}\cutoffint_{T_x^*U} d\xi\, \frac{ D_\xi^\alpha\sigma(x,\xi+\theta \eta)}{\alpha!} \eta^\alpha.
$$
Since  $\sigma$ has non integer order symbol,  neither  has $D^\alpha\sigma$ an integer order. After integrating over $U$,  the terms corresponding to $\vert \alpha\vert \neq 0$ vanish by the integration by parts formula, as a result of which we are left with the 
 $\vert \alpha\vert=0$ term and 
$$
\cutoffint_{T^*U} \sigma(x,\xi+\eta)\, dx\, d\xi 
= \cutoffint_{T^*U}  \sigma(x,\xi)\, dx\, d\xi.$$\endsquare
\subsection{The cosphere cochain as a $B_0$-coboundary }
\begin{defn}
 Let the {\rm cosphere} $k$-cochain denote the $(k+1)$-linear form on $CS_{\rm com}(U)$ 
$$
\psi_{k}(\sigma_0, \cdots, \sigma_{k})= \int_{S^*U}j^*\left( \sigma_0\ast d\, \sigma_1\wedge_\ast \cdots \wedge_\ast d\,\sigma_{k}\right)_0
$$
for all $ \sigma_0, \cdots, \sigma_{k} \in CS_{\rm com}(U)$.
\end{defn}
Since $\psi_{k}$ vanishes for $k<2n-1$, we shall focus on $\psi_{2n-1}$. \\
We introduce a cochain on $CS_{com}(U)$ built from cut-off integrals of classical symbol valued forms:
\begin{defn}For any  $\sigma_0, \cdots, \sigma_k \in CS_{\rm com}(U)$ we set
$$\chi_{k}^{\rm cut-off}(\sigma_0, \cdots, \sigma_{k})= \cutoffint_{T^*U}\left(\sigma_0\ast d\, \sigma_1\wedge_\ast \cdots \wedge_\ast d\,\sigma_{k}\right)_0.$$
\end{defn}
\begin{rk}
$\chi_k^{\rm cut-off}$ vanishes for $k<2n$ so that we focus on the $2n$-cochain  $\chi_{2n}^{\rm cut-off}$.
\end{rk} 
By Stokes' formula for cut-off integrals on  non integer order symbol valued
forms, we have (with $B_0$ as in Appendix A):
\begin{eqnarray*}
B_0\chi_{2n}^{\rm cut-off}(\sigma_0, \cdots, \sigma_{2n-1})&=& \chi_{2n}^{\rm cut-off}(1,\sigma_0, \cdots, \sigma_{2n-1 })\\
&=& \cutoffint_{T^* U} \left(d\sigma_0\wedge_\ast \cdots\wedge_\ast d\sigma_{2n-1}\right)_0\\
&=& \cutoffint_{T^* U} d\left(\sigma_0\wedge_\ast d\sigma_1\wedge_\ast \cdots\wedge_\ast d\sigma_{2n-1}\right)_0\\
&=&0
\end{eqnarray*}
whenever the sum of the orders of the $\sigma_i$'s is non integer.\\
 However, 
$\chi_{2n}^{\rm cut-off}$ is  not cyclic in general; the obstruction
to its cyclicity is measured by the cosphere cochain.  
\begin{prop}\label{prop:Bcutoff}
\begin{eqnarray*}
B_0\chi_{2n}^{\rm cut-off}(\sigma_0, \cdots, \sigma_{2n-1})&=& \chi_{2n}^{\rm cut-off}(1,\sigma_0, \cdots, \sigma_{2n-1 })\\
&=& \psi_{2n-1}(\sigma_0, \cdots, \sigma_{2n-1})
\end{eqnarray*}
for any $\sigma_0, \cdots, \sigma_{2n-1}\in CS_{com}(U).$\\
It vanishes whenever the $\sigma_i$'s have orders which sum up to a non integer. 
\end{prop}
{\bf Proof:}
\begin{eqnarray*}
B_0\chi_{2n}^{\rm cut-off}(\sigma_0, \cdots, \sigma_{2n-1})&=& \chi_{2n}^{\rm cut-off}(1,\sigma_0, \cdots, \sigma_{2n-1 })\\
&=& \cutoffint_{T^*U} \left(d\, \sigma_0\wedge_\ast \cdots \wedge_\ast
 d\,\sigma_{2n-1}\right)_0\\
&=& {\rm fp}_{R\to\infty}\int_{B^*(0, R)} d\,\left( \sigma_0 \ast d\, \sigma_1  \cdots \wedge_\ast d\,\sigma_{2n-1}\right)_0\\
&=& {\rm fp}_{R\to\infty}\int_{S^*(0, R)} \left( \sigma_0\ast d\sigma_1 \cdots \wedge_\ast d\,\sigma_{2n-1}\right)_0\\
&=&\int_{S^*(0, 1)} \left( \sigma_0\ast d\sigma_1\wedge_\ast\cdots \wedge_\ast d\,\sigma_{2n-1}\right)_0\\
&=&\psi_{2n-1}\left( \sigma_0,\sigma_1, \cdots ,\sigma_{2n-1}\right).\\
\end{eqnarray*}
\endsquare
\section{The Wodzicki residue extended to forms as a complex residue}
We first recall how the ordinary residue density on symbols can be interpreted as a complex residue via cut-off integrals of symbols. 
\subsection{The Wodzicki residue density on symbols  as a complex residue}
Recall that given an open subset $U\subset \R^n$ (resp. an $n$-dimensional manifold $M$),  for any real number $m$ the class $CS_{com}^m(U)$  of classical symbols of order $m$ with compact support on $U$ (resp.  of classical symbols of order $m$)  can be equipped with a natural Fr\'echet topology  so that 
$ \bigcup_{m\in \R} CS_{com}^m(U)$   comes equipped with an inductive limit Fr\'echet topology. We first recall  the notion of holomorphic regularisation (see e.g.\cite{P} for a review of various regularisations):
\begin{defn}
A  holomorphic regularisation procedure on $CS_{com}(U)$  
is a  map 
\begin{eqnarray*}
{\cal R}: CS_{com}(U)&\to & {\rm Hol}\,\left( CS_{com}(U)\right)\\
\sigma  &\mapsto & \sigma(z)
\end{eqnarray*}
where ${\rm Hol}\left( CS_{com}(U)\right)$ is the algebra of holomorphic maps with values in $CS_{com}(U)$, 
such that 
\begin{enumerate}
\item  $\sigma(0)=\sigma$, 
\item  $\sigma(z)$ has holomorphic  order  $\alpha(z)$ (in particular,
  $\alpha(0)$ is equal to the order of $\sigma$)  such that
  $\alpha^\prime(0)\neq 0$.
\end{enumerate}
By holomorphic map we mean that each 
  positively homogeneous component $\sigma_{\alpha(z)-j}(z)$ is holomorphic and that 
 for any integer $N\geq 1$ the remainder
$$\sigma_{(N)}(z)(x, \xi):= \sigma(z)(x,\xi)- \sum_{j=0}^{N-1}
\sigma_{\alpha(z)-j}(z)(x, \xi)$$ is holomorphic in $z$ as an
element of $\Ci( U \times \R^n)$ with $k^{{\rm
th}}$ $z$-derivative
\begin{equation}\label{e:kthderivlogclassical}
\sigma^{(k)}_{(N)}(z)(x, \xi) := \partial_z^k(\sigma_{(N)}(z)(x, \xi))
\in {\rm S}^{\alpha(z)-N + \e}(U,V)
\end{equation}
for any $\e>0$.
\end{defn}
A first example of holomorphic regularisation is  the well known Riesz regularisation, which sends  a
classical symbol $\sigma$ of order $m$ to $$ \sigma(z)(x, \xi):= \sum_{j=0}^N \chi(\xi) \, \sigma_{m-j}(x, \xi)\cdot \vert
\xi \vert^{-z}+\sigma_{(N)}(x, \xi)$$ with the notations of 
(\ref{eq:classicalsymb})
 and where  $N$ is chosen large enough so that
$m-N<-n$.  Generalisations  of the type  $$\sigma\mapsto \sigma(z)(x,
\xi):= H(z) \sum_{j=0}^N \chi(\xi) \, \sigma_{\alpha-j}(x, \xi)\cdot \vert
\xi \vert^{-z}+\sigma_{(N)}(x, \xi)$$  where $H$ is a
holomorphic function such that $H(0)=1$   include dimensional
regularisation which arises in physics (see \cite{P}).
\begin{rk}
\end{rk}
\begin{prop}\label{prop:KV} \cite{G}, \cite{KV}, \cite{L}
Given a holomorphic regularisation procedure ${\cal R}:\sigma \mapsto \sigma(z)$ on $CS_{com}(U)$  and any symbol $\sigma\in CS_{com}(U)$, for any $x\in U$ the map $z\mapsto \cutoffint_{T_x^*U}\, d\xi \,  \sigma(z)$ (resp.  $z\mapsto \cutoffint_{T^*U} \, dx\, d\xi\, \sigma(z)$) is meromorphic with simple poles at points in $\alpha^{-1}([-n, +\infty[\, \cap\, \Z)$ where $\alpha$ is the order of $\sigma(z)$. Moreover for any $x\in U$
$${\rm Res}_{z=0} \cutoffint_{T_x^*U} \sigma(z)(x,\xi) \, d\xi= -\frac{1}{\alpha^\prime(0)} {\rm res}_x(\sigma(0)),$$
respectively 
$${\rm Res}_{z=0} \cutoffint_{T^*U} \sigma(z)(x,\xi) \, d\xi= -\frac{1}{\alpha^\prime(0)} {\rm res}(\sigma(0)).$$
\end{prop}
On the grounds of this proposition we set:
\begin{defn}
The ${\cal R}$-regularised integral of $\sigma\in CS_{com}(U)$    is defined by:
\begin{eqnarray*}
\int_{T^*U}^{{\cal R}} \sigma(x, \xi) \, d\xi &:=& {\rm fp}_{z=0} \cutoffint_{T^*U} \sigma(z)(x, \xi)\, d\xi\\
&:=& \lim_{z\to 0}\left(  \cutoffint_{T^*U}d\xi \,  \sigma(z)(x, \xi)- \frac{1}{z} {\rm Res}_{z=0}\cutoffint_{T^*U} d\xi\, \sigma(z)(x, \xi)\right)
\end{eqnarray*}
\end{defn}
{\bf Proof of the proposition:} We identify $T_x^*U$ with $\R^n$ using a coordinate chart.  From equation (\ref{eq:finitepart}) we have 
\begin{eqnarray*}
\cutoffint_{\R^n} \sigma(z)(x,\xi) \, d\xi&= & \int_{B(0, 1)}
 \sigma(z)(x, \xi) \, d\xi\\
&\hskip -15mm - &\hskip -15mm \sum_{i=0, \, \alpha(z)-i+n\neq 0 }^{N} \, \frac{1}{\alpha(z)-i+n} \, \int_{S(0, 1)}\sigma_{\alpha(z)-i}(z)(x, \xi) \, d\xi\\
&\hskip -15mm +&\hskip -15mm   \int_{\R^n}
 \sigma_{(N)}(z)(\xi) \, d\xi\\
&=& \int_{B(0, 1)} \sigma(z)(x,\xi) \, d\xi\\
&\hskip -15mm - &\hskip -15mm \sum_{i=0, \, \alpha(z)-i+n\neq 0 }^{N} \,\frac{1}{\alpha(0)-i+n +\alpha^\prime(0) \,z +o(z)} \, \int_{S(0, 1)}\sigma_{\alpha(z)-i}(z)(x,\xi) \, d\xi\\
 &\hskip -15mm + &\hskip -15mm  \int_{\R^n}
 \sigma_{(N)}(z)(x,\xi) \, d\xi,\\
\end{eqnarray*}
where we have written  $\alpha(z)=\alpha(0)+\alpha^\prime(0).z +o(z)$. 
As a consequence, we have that:
\begin{eqnarray*}
{\rm Res}_{z=0} &&\hskip -8mm\cutoffint_{\R^n} \sigma(z)(x,\xi) \, d\xi\\
&= &{\rm Res}_{z=0}\int_{B(0, 1)}
 \sigma(z)(x,\xi) \, d\xi\nonumber\\
&-& {\rm Res}_{z=0}\sum_{i=0 }^{K_N}\, \frac{1}{\alpha(0)-i+n +\alpha^\prime(0) \,z +o(z)} \, \int_{S(0, 1)}\sigma_{\alpha(z)-i}(z)(x,\xi) \, d\xi\\
&+&  {\rm Res}_{z=0} \int_{\R^n}
 \sigma_{(N)}(z)(x,\xi) \, d\xi\\
&= &- \frac{1}{\alpha^\prime(0)} \, \int_{S(0, 1)}\sigma_{-n}(0)(x,\xi) \, d\xi\\
&=& - \frac{1}{\alpha^\prime(0)} \,{\rm res}_x(\sigma(0)).\\
\end{eqnarray*}
\endsquare
This result  extends to classical symbol valued forms.
\subsection{Cut-off integrals of holomorphic families of symbol valued forms}
\begin{defn}
A  holomorphic regularisation procedure on $\Omega\, CS(U)$  
is a  map 
\begin{eqnarray*}
{\cal R}: \Omega\, C_{com}(U)&\to &\Omega\,  {\rm Hol}\,\left( CS_{com}(U)\right)\\
\omega  &\mapsto & \omega(z)
\end{eqnarray*}
where 
\goodbreak
\begin{eqnarray*}
\Omega\, {\rm Hol}\left( CS_{com}(U)\right)&:=&\{z\mapsto \omega (z)= \sum_{I, J} \omega_{IJ} dx_I\wedge d\xi_J \in \Omega \, CS_{com}(U),\\
& \quad & z\mapsto \omega_{IJ} (z) \quad {\rm lies} \quad {\rm in} \quad {\rm Hol}\,CS(U)\\
&{}& \quad {\rm for } \quad {\rm all} \quad {\rm multi-indices} 
\quad I, J\}\\
\end{eqnarray*}
and 
\begin{enumerate}
\item  $\omega(0)=\omega$, 
\item  $\omega(z)$ has holomorphic  order  $\alpha(z)$ (in particular, $\alpha(0)$ is equal to the order of $\omega$)  such that $\alpha^\prime(0)\neq 0$.\\
\end{enumerate}
\end{defn}
\begin{rk} 
Clearly, any holomorphic regularisation ${\cal R} $ on $CS_{com}(U)$ induces one on $\Omega\, CS_{com}(U)$ setting:
$${\cal R}(\omega) = \sum_{I, J} {\cal R}(\omega_{IJ}) dx_I\wedge d\xi_J.$$
\end{rk}
\begin{thm}\label{thm:KVforms} 
Given a holomorphic regularisation procedure ${\cal R}:\omega \mapsto \omega(z)$ on $\Omega \, CS_{com}(U)$  induced by a regularisation ${\cal R}: \sigma\mapsto \sigma(z)$ on $CS_{com}(U)$   and any symbol valued form $\omega\in \Omega\, CS_{com}(U)$,  the map $z\mapsto \cutoffint_{T_x^*U}\,  \omega(z)$ (resp.  $z\mapsto \cutoffint_{T^*U} \,  \omega(z)$) is meromorphic with simple poles at points in $\alpha^{-1}([-n, +\infty[\, \cap\, \Z)$ where $\alpha$ is the order of $\omega(z)$. Moreover for any $x\in U$
$${\rm Res}_{z=0} \cutoffint_{T_x^*U} \omega(z)(x,\xi) = -\frac{1}{\alpha^\prime(0)} {\rm res}_x(\omega(0)),$$
respectively 
$${\rm Res}_{z=0} \cutoffint_{T^*U} \omega(z)(x,\xi) = -\frac{1}{\alpha^\prime(0)} {\rm res}(\omega(0)).$$
\end{thm}
On the grounds of this theorem,  we set the following definition:
\begin{defn} The ${\cal R}$-regularised integral of  $\omega\in \Omega CS_{com}(U)$ is defined by:
\begin{eqnarray*}
\int_{T^*U}^{{\cal R}} \omega(x, \xi)  &:=& {\rm fp}_{z=0} \cutoffint_{T^*U} \omega(z)(x, \xi)\\
&:=& \lim_{z\to 0} \left(  \cutoffint_{T^*U} \omega(z)(x, \xi)- \frac{1}{z} {\rm Res}_{z=0}\cutoffint_{T^*U} \omega(z)(x, \xi)\right). 
\end{eqnarray*}
\end{defn}
{\bf Proof of the theorem:} The result  follows from applying Proposition 
\ref{prop:KV} to each component $\omega_{IJ}(z)$ of the form 
$\omega(z)=\sum_{IJ} \omega_{IJ} (z) dx_I\wedge d\xi_J$. The symbol valued form $\omega_{IJ}(z)$ has order $\alpha_{IJ}(z)=\alpha(z)-\vert J\vert$ so that $\alpha_{IJ}^\prime(0)= \alpha^\prime(0)$. Since $z\mapsto \cutoffint_{T_x^*U} \omega_{IJ}(z)$ is meromorphic with simple poles so is $z\mapsto \cutoffint_{T_x^*U} \omega(z)$ and we have 
\begin{eqnarray*}
{\rm Res}_{z=0} \cutoffint_{T_x^*U} \omega(z)(x,\xi)&=& \sum_{IJ} 
{\rm Res}_{z=0} \cutoffint_{T_x^*U} \omega_{IJ}(z)(x,\xi)\, dx_I\wedge d\xi_J \\
& =& -\sum_{IJ}\frac{1}{\alpha_{IJ}^\prime(0)} {\rm res}_x(\omega_{IJ}(0))\, dx_I\wedge d\xi_J\\
&{}& {\rm by} \quad {\rm Proposition} \quad \ref{prop:KV}\\
&=& -\frac{1}{\alpha^\prime(0)}\sum_{IJ}{\rm res}_x(\omega_{IJ}(0))\, dx_I\wedge d\xi_J\\
& =& -\frac{1}{\alpha^\prime(0)}{\rm res}_x(\omega(0)),
\end{eqnarray*}
\endsquare
Stokes' formula holds as an equality of meromorphic functions:  
\begin{thm}\label{thm:Stokesreg}
Given a holomorphic regularisation procedure ${\cal R}:\omega \mapsto \omega(z)$ on $\Omega \, CS_{com}(U)$  induced by a regularisation ${\cal R}: \sigma\mapsto \sigma(z)$ on $CS_{com}(U)$   and any symbol valued form $\omega\in \Omega\, CS_{com}(U)$, we have the following equality of meromorphic functions:
$$\cutoffint_{T_x^*U}d\, \left( \omega(z)\right)=0.$$
\end{thm}
{\bf Proof:}
Since $\omega(z)$ has non integer order outside a discrete set of complex numbers, and since by Theorem \ref{thm:Stokescutoff},
Stokes' property holds for non integer order symbols valued forms, the statement holds outside this discrete set of poles. The meromorphicity of the  function $z\mapsto \cutoffint_{T_x^*U}d\, \left( \omega(z)\right)$ proved in Theorem \ref{thm:KVforms} then yields the expected equality of meromorphic functions.
\begin{rk}
\begin{itemize}
\item  This statement in the case of dimensional regularisation and  transposed to forms built from Feynman type functions as in \cite{E}   corresponds  to Proposition 12 of \cite{E}. 
\item In general, $$\cutoffint_{T_x^*U}^{\cal R } d\,  \omega\neq 0$$ since  exterior differentiation and  regularisation ${\cal R}$ do not ``commute''. 
\end{itemize}
\end{rk}
\section*{Appendix A} We recall here a few
definitions borrowed from non commutative geometry see e.g.\cite{C},
\cite{GVF}. Let  $({\cal A},\star)$  be an associative algebra over some ring
$R$ with unit $1$.  The  space  $C^n\left( {\cal A}, R\right)$  of
$R$-valued  $n+1$-linear forms on ${\cal A}$ corresponds to the
space of $n$-cochains on ${\cal A}$. Equivalently, these spaces can be seen as
 spaces of  $R$-multilinear  $n$-forms on ${\cal A}$ with values in the
$R$-algebraic  dual ${\cal A}^*$, seen as an ${\cal A}$-bimodule, where for
$\chi \in {\cal A}^*$ we put $a^\prime \chi(a) a^{\prime \prime}=
\chi(a^{\prime\prime} a a^{\prime})$.\\  Following \cite{C} we  define the
operators $B_0$ and $B$ acting on cochains: 
\begin{defn} Let 
\begin{eqnarray*}  B_0:C^n ({\cal A})&\to & C^{n-1}({\cal A})\\ \chi&\mapsto &
B_0\chi(a_0, \cdots, a_{n-1}):=\chi(1, a_0, \cdots, a_{n-1})- (-1)^n  \chi (
a_0, \cdots,  a_{n-1},1). \end{eqnarray*} Let   $B:= {\cal A}\, B_0$  where  
${\cal A}$ denotes cyclic antisymmetrisation in all variables  so that 
$$B\chi(a_0, \cdots, a_{n-1})=  \sum_{i=0}^{n-1} (-1)^i \chi(1,a_i,a_{i+1},
\cdots)- (-1)^n \sum_{i=0}^{n-1} (-1)^i \chi(a_i,a_{i+1}, \cdots, a_{i-1},
1)$$ \end{defn} One can check that $B^2=0$ so that $B$ defines a homology on
$C^{\bullet} ({\cal A})$ \cite{C}.  \begin{defn} The Hochschild coboundary for
the product $\star$  of an  $n$-cochain $\chi$ is defined by:
$$b_\star\chi(a_0,\cdots, a_{n+1})= \sum_{j=0}^n (-1)^j \chi(a_0, \cdots, a_j
\star a_{j+1}, \cdots, a_{n+1})+ (-1)^{n+1} \chi(a_{n+1} \star a_0, \cdots,
a_n).$$ \end{defn} It satisfies the condition $b_\star^2=0$ and hence defines
a cohomology called the {\it Hochschild cohomology} of $({\cal A}, \star)$. 
\begin{defn}  An $n$-dimensional  cycle  is  given by a triple $\left(\Omega,
d, \int\right)$ where  $\Omega$ is a graded differential algebra on $\C$ 
equipped with the  differential $d$ such that $d^2=0$ and  $\int: \Omega^n\to
\C$  is a {\it closed} graded trace i.e.  $\int$ is  a linear map  which, when
extended  to $\Omega$ by $0$, satisfies  $$\int \alpha\wedge \beta =
(-1)^{\vert \alpha\vert \cdot \vert \beta\vert}\cdot \int \beta \wedge \alpha,
\quad \int d\beta=0 \quad\forall \, \beta \in \Omega^{n-1}({\cal A}).$$  An
$n$-cycle on  an algebra ${\cal A}$ on $\C$ is a cycle $(\Omega, d, \int)$ 
together with a homomorphism $\rho: {\cal A}\to \Omega^0$.  The character
$\chi_n$  of an $n$-cycle is defined by: $$\chi_n(a_0, \cdots, a_n)= \int
\rho(a_0)\, d\rho(a_1)\, \cdots \, d\rho(a_{n})\quad \forall a_i\in {\cal
A}.$$  \end{defn} Let us also recall that the character of a cycle  has the
following properties: \begin{enumerate} \item $\chi_n$ is cyclic i.e. 
$$\chi_n(a_0, \cdots, a_n)= (-1)^{n} \chi_n(a_1, \cdots, a_n, a_0),\quad
\forall a_i\in {\cal A}$$ \item $\chi_n(1, a_1, \cdots, a_n)=0\quad \forall
a_i\in {\cal A}$.  \item 
$b\, \chi_n=0$ where $b$ is the Hochschild coboundary associated with the
product on ${\cal A}$.  \end{enumerate}

 \bibliographystyle{plain}  \end{document}